```
\documentclass[12pt]{article}
\usepackage{amsmath, amsfonts, amsthm, amssymb,epsfig}
%\input{macros.TEX}

\begin{document}

\title{Pythagorean Boxes with Primitive Faces}

\author{Konstantine D. Zelator\\
P.O. Box 4280\\
Pittsburgh, PA  15203}

\maketitle

\section{Introduction} In their paper ``Pythagorean Boxes", Raymond A.
Beauregard and E. R. Suryanarayan define the concept or notion of a
Pythagorean rectangle as one with integer sides and integer diagonals
(see \cite{r1}); they also introduce the concept of a Pythagorean Box as
a rectangular three-dimensional box whose edges and inner diagonals are
integers (or more precisely, they have integer lengths).  As in that
paper, the abbreviation {\bf PB} will simply stand for ``Pythagorean
Box"; also in this article, the abbreviated notion {\bf PR} will stand
for ``Pythagorean Rectangle".

In the Beauregard and Suryanarayan paper, it was shown that there exist
infinitely many {\bf PB's} with a square base and height equal to 1.

In this paper, we present a method and formulas that generate infinitely
many PB's that contain a pair of opposite (and hence congruent) PR's
which are primitive; a Pythagorean Rectangle is primitive if the  four
congruent pythagorean triangles contained therein are primitive.

In Figure 1 a PB is illustrated; it has edge lengths the integers
$x,y,z$; $t$ is the integer length of the four inner diagonals and
$s,u,w$ are the lengths of the twelve face diagonals (four of them having
length $s$, four having length $u$, and four having length $w$).
\newpage

\ \ \
\begin{center}
\epsfig{file=pythboxfig1.eps,width=2.5in}
\end{center}

\vspace{.15in}

\noindent In any PB the four positive integers $x,y,z,t$ must satisfy the
equation

\begin{equation}
x^2+y^2+z^2 = t^2 \label{E1}
\end{equation}

If a PB contains a face which is a PR, then of course the oppositive face
must also be a PR; there are three pairs of congruent faces.  Suppose
```

that the two opposite faces characterized by the triple $(x,z,w)$ are PR's. That means that in addition to $x$ and $z$, $w$ must also be an integer. In other words, the triple $(x,z,w)$ must be a Pythagorean one; moreover, if $(x,z)=1$, then $(x,z,w)$ will be a primitive Pythagorean triple. Also, the relations $x^2+z^2=w^2$ and $(x,z)=1$ are equivalent to $x^2+z^2=w^2$ and $(x,z)=1=(x,w) = (w,z) = (x,z,w)$.

In terms of notation, in what is to follow, the PB at hand will be characterized by the septuple $(x,y,z,t,s,u,w)$ while the three pairs of opposite faces will be characterized by the three triples $(x,z,w),\ (x,y,s)$ and $(y,z,u)$; in accordance with Figure 1.

There are three results in this paper. In Result 1, we derive certain explicit conditions that a PB must satisfy, if it possesses two pairs of (opposite) primitive PR's. In Result 2, we show that if similar conditions are satisfied then infinitely many PB's containing a pair of (opposite) PR's can be generated. In Result 3, we prove that there are no PB's with a square base and with a face which is a primitive PR.

\section{Listing of Results}

\noindent {\bf Result 1:} If a PB, characterized by the septuple $(x,y,z,t,s,u,w)$, has two pairs of faces, characterized by the triples $(x,z,w)$ and $(y,z,u)$, and which are primitive PR's, then there exist postive integers $M_1,N_1,M_2,N_2,\delta$ such that

\framebox{
\parbox[t]{4.75in}{\begin{center}$M^2_1 +N^2_1+\delta \cdot N^2_2 = \delta \cdot M^2_2$\end{center}
\begin{center}$(M_1,N_1) = 1 = (M_2,N_2),\ M_1+N_1 \equiv M_2+N_2 \equiv 1({\rm mod}\ 2),$\end{center}
\begin{center}$x=2M_1N_1,\ y = 2\delta M_2N_2,\ z=M^2_1-N^2_1,\ w=M^2_1 +N^2_1$,\end{center}
\begin{center}$u = \sqrt{(M^2_1-N_1)^2+(2\delta M_2N_2)^2};\ u$\end{center}
being a positive integer, $t=\delta (M^2_2+N^2_2)$, and $\delta$ an odd positive integer.}}

\vspace{.15in}

\noindent{\bf Result 2:} Suppose that the positive integers $M_1,N_1,M_2,N_2$ satisfy the conditions, $(M_1,N_1)=1,\ M_1+N_1\equiv 1({\rm mod}\ 2),\ M_1 > N_1$ and $M^2_1 +N^2_1+N^2_2 = M^2_2$. Define the integers, $x = 2M_1N_1,\ y=2M_2N_2,\ z=M^2_1-N^2_1,\ w=M^2_1+N^2_1$, and $t = M^2_1 +N^2_1 +2N^2_2$; and the real numbers, $s = \sqrt{x^2+y^2} = \sqrt{(2M_1N_1)^2+(2M_2N_2)^2}$ and $u = \sqrt{y^2+z^2} = \sqrt{(2M_2N_2)^2+(M^2_1-N^2_1)^2}$. Under the above assumptions, the 3 dimensional rectangular parallelepiped characterized by the septuple $(x,y,z,t,s,u,w)$ is a PB with a pair of opposite faces which are primitive PR's; those characterized by the triple $(x,z,w)$.

\vspace{.15in}

\noindent{\bf Result 3:}  There exist no PB's with a square base and which possess a face which is a primitive PR.  In other words, there do not exist any PB's of the form $(x,x,z,t,s,u,w)$ (i.e., $(x,y,z,t,s,u,w)$ with $x=y$) and with a face which is a primitive PR.  Equivalently, if such a PB has a face which is a PR, then that face must be nonprimitive.

With respect to Result 2, note that if $u$ is an integer, then the two (opposite) faces corresponding to the triple $(y,z,u)$ will be PR's; if in addition the condition $(2M_2N_2,M^2_1-N^2_1)=1$ (which is equivalent to $(M_2N_2,M^2_1-N^2_1)=1$ since $M^2_1-N^2_1$ is an odd integer) is satisfied, the two faces characterized by the triple $(y,z,u)$, will be in fact primitive PR's.

Let us make a few remarks on the other two results listed above.  In the first result, note that for $\delta = 1$, the numbers $M_1,N_1,M_2$, and $N_2$ form a solution to equation (\ref{E1}); specifically the quadruple $(M_1,N_1,N_2,M_2)$ forming a solution to equation (\ref{E1}).  If we use Figure 1 as our reference, then if the illustrated PB has lateral faces which are all PR's, then those formulas describing the integer lengths $x,y,z,t,u,$ and $w$ (in terms of the four integers $M_1,N_1,M_2,N_2$) must hold true; also note that, since $u$ is required to be an integer, the radicand quantity in the formula for $u$ must be a perfect square; in other words, the integer $(M^2_1-N^2_1)^2+(2M_2N_2)^2$ must equal an integer square.  One of the interesting features of these formulas, is that the quadruple $(M_1,N_1,N_2,M_2)$ is a solution to an equation similar to (\ref{E1}).  The general solution to the diophantine equation (\ref{E1}) is well-known and can be found in \cite{r2}. Actually, we make use of the general solution to equation (\ref{E1}) in the proof of Result 1.

In the proof of Result 3, we make use of two results from number theory. The first result we employ in the proof of Result 3 is the general solution, in positive integers, to the diophantine equation $X^2+2Y^2=Z^2$; this can be found in \cite{r3} and is well-known.  The second result we use is the fact that there exist no two positive integers both of whose sum and difference of their squares, are also integer squares; it can be found in W. Sierpinski's ``Elementary Theory of Numbers" (refer to \cite{r2}).  However, this is perhaps not as widely known result, we offer it as Lemma 1 and its proof, just after the proofs (of the three results) section.

\vspace{.15in}

\noindent{\bf Proofs}

\vspace{.15in}

\noindent {\bf 1.} {\it Proof of Result 1}

\vspace{.15in}

Since $(x,y,w)$ and $(y,z,u)$ are primitive Pythagorean triples, the positive integers $x,y,z,t,u$, and $w$ must satisfy the conditions,
\setcounter{equation}{0}

```
\begin{equation} x^2+y^2+z^2 = t^2 \label{E1a}
\end{equation}

\begin{equation}
x^2+z^2=w^2 \label{E2}
\end{equation}

\begin{equation}
y^2+z^2 = u^2 \label{E3}
\end{equation}

\begin{equation}
(x,z)=1=(y,z) \label{E4}
\end{equation}


Just from equation (\ref{E1a}) alone, an easy argument modulo 4 shows
that at least two of the three integers $x,y,z$ must be even; this
conclusion combined with condition (4) implies that $x$ and $y$ must be
even, while $z$ must be odd; and consequently all three integers $t,w$
and $u$ must also be odd.  We now employ the general solution to equation
(\ref{E1}) which can be found in reference \cite{r2}.

\vspace{.15in}

\noindent Equations (5) \framebox{
\parbox[c]{3.5in}{$x=2\ell, \ y = 2m,\ z = {\displaystyle
\frac{\ell^2+m^2-n^2}{n}},\ t = {\displaystyle
\frac{\ell^2+m^2+n^2}{n}}$, for positive integers $\ell, m,n$; $n$ being
a divisor of $\ell^2+m^2$.}}

\vspace{.15in}

Since $(x,z,w)$ is a primitive Pythagorean triple we must have,

\vspace{.15in}

\noindent Equations (6): \framebox{
\parbox[c]{3.5in}{$x=2M_1N_1,\ z=M^2_1-N^2_1,\ w = M^2_1+N^2_1$, for
positive integers $M_1,N_1$ such that $(M_1,N_1)=1,\ M_1>N_1$, and
$M_1+N_1 \equiv 1({\rm mod}\ 2)$.}}

\vspace{.15in}

Likewise, since $(y,z,u)$ is a primitive Pythagorean triple,

\vspace{.15in}

\noindent Equations (7): \framebox{
\parbox[c]{3.5in}{$y = 2M_0,N_0,\ z=M^2_0-N^2_0,\ u = M^2_0 +N^2_0$, for
positive integers $M_0,N_0$ such that $(M_0,N_0)=1,\ M_0 > N_0$, and $M_0
+ N_0 \equiv 1({\rm mod}\ 2)$.}}
```

```latex
\vspace{.15in}

Combining equations (5), (6), and (7) we see that,

$$
\ell = M_1N_1,\ m = M_0,N_0,\ M^2_1 - N^2_1 = M^2_0 - N^2_0 =
{\displaystyle \frac{M^2_1N^2_1+M^2_0N^2_0 - n^2}{n}}.
$$

In particular, we obtain,

\setcounter{equation}{7}

\begin{equation}
\begin{array}{ll}
n(M^2_1-N^2_1) = M^2_1N^2_1 +M^2_0N^2_0 - n^2;\\
\\
n^2+n(M^2_1-N^2_1) - (M^2_1N^2_1 +M^2_0N^2_0)=0
\end{array}\label{E8}
\end{equation}

Since the integer $n$, according to (\ref{E8}), is a root to the
quadratic equation with integer coefficients, $X^2+(M^2_1-N^2_1)X -
(M^2_1N^2_1+M^2_0N^2_0)=0$, it follows that the discriminant of the
latter equation must be a perfect or integer square:

\begin{equation}
\begin{array}{ll}
(M^2_1-N^2_1)^2 + 4(M^2_1N^2_1 + M^2_0N^2_0) = K^2;\\
\\
(M^2_1+N^2_1)^2 + (2M_0N_0)^2 = K^2.
\end{array}
\label{E9}
\end{equation}

Equation (\ref{E9}) shows that the triple $(M^2_1+N^2_1, 2M_0N_0,K)$ is a
Pythagorean one, though not necessarily primitive.  Also note that since
$w=M^2_1+N^2_1$ (from equations (6)) and $y = 2M_0N_0$ (from equations
(7)), and in view of equation (\ref{E1}), (2), and (3), the positive
integer $K$ is none other than $t$; $K=t$.  We have (since $M^2_1+N^2_1$
is odd and $2M_0N_0$ is even),

\vspace{.15in}

\framebox{
\parbox[c]{3.5in}{\begin{center}
$M^2_1+N^2_1 = \delta (M^2_2 - N^2_2)$
\end{center}
\begin{center}
$2M_0N_0 = \delta \cdot (2M_2N_2)$
\end{center}
\begin{center}
$K=t=\delta \cdot (M^2_2+N^2_2),$
```

```
\end{center}
where $\delta$ is an odd integer and the positive integers $M_2,N_2$
satisfy $M_2>N_2,\ (M_2,N_2)=1$ and $M_2+N_2 \equiv 1({\rm mod}\ 2)$.}}
\hfill Equations (10)

\vspace{.15in}

\noindent From (6) we know that $x = 2M_1N_1, \ z=M^2_1-N^2_1,\ w =
M^2_1+N^2_1$, with $M_1+N_1 \equiv 1({\rm mod}\ 2),\ M_1>N_1$, and
$(M_1,N_1)=1$.  Moreover, from (7) and (10) we deduce that $y=\delta
\cdot (2M_2N_2)$, with $\delta$ being an odd integer and the positive
integers $M_2$ and $N_2$ satisfying $M_2>N_2,\ (M_2,N_2)=1$, and $M_2+N_2
\equiv 1({\rm mod}\ 2)$.  And from (3) $\Rightarrow u = \sqrt{y^2+z^2}$;
and thus, $y=\sqrt{\left(M^2_1-N^2_1\right)^2+\left(2\delta
M_2N_2\right)^2}$. Finally, the first equation in (10) implies
$M^2_1+N^2_1+\delta N^2_2 = \delta M^2_2$. The proof is complete. \hfill
\rule{2mm}{2mm}

\section{Proof of Result 2}

Obviously, $\left(M^2_1-N^2_1\right)^2+\left(2M_1N_1\right)^2 =
\left(M^2_1+N^2_1\right)^2$ and thus by definition, $x^2+z^2 = w^2$; and
since $\left(M^2_1-N^2_1, 2M_1N_1\right) =1$, on account of the
conditions $\left(M_1,N_1\right) = 1$ and $M_1+N_1 \equiv 1({\rm mod}\
2)$, we conclude that $(x,z)=1$ and thus $(x,z,w)$ is a primitive
pythagorean triple and therefore the two opposite rectangulat faces which
correspond to or are characterized by the triple $(x,z,w)$, are actually
primitive PR's.  To finish the proof, we must show that actually the
septuple $(x,y,z,t,s,u,w)$ truly corresponds to a PB.

We must establish the fundamental condition,

\setcounter{equation}{10}
\begin{equation}
x^2+y^2+z^2=t^2 \label{E11}
\end{equation}

\noindent This we do by establishing an equivalent true statement, by way
of direct computation:

$$
\begin{array}{rcl}
(1) & \Leftrightarrow  & \left(2M_1N_1\right)^2 + \left( 2M_2N_2\right)^2
+ \left(M^2_1-N^2_1\right)^2 = \left(M^2_1 +N^2_1 +2N_2\right)^2 \\
\\
& \Leftrightarrow & 4M^2_1N^2_1 +4M^2_2N^2_2 +M^4_1 - 2M^2_1N^2_1 +N^4_1
= M^4_1+N^4_1 \\
\\
& & +4N^4_2+2M^2_1N^2_1 +4M^2_1N^2_2 + 4N^2_1N^2_2\\
\\
& \Leftrightarrow & M^2_2N^2_2= N^4_2 +M^2_1 N^2_2 +N^2_1 N^2_2\\
\\
&\Leftrightarrow&({\rm since}\ N_2\neq 0)\ M^2_2 = N^2_2 +M^2_1 +N^2_1,
\end{array}
```

```
$$

\noindent which is true by hypothesis.  \hfill \rule{2mm}{2mm}

\section{Proof of Result 3}

We will show that the assumption of the existence of a PB of the form
$(x,x,z,t,s,u,w)$, in which a face is a primitive PR, leads to a
contrdiction.  There are three pairs of opposite faces characterized by
the triples $(x,x,s),\ (x,z,w)$, and $(x,z,u)$.   Obviously $u = w$, so
there are four congruent lateral faces.  If a face were a primitive PR,
so would be its opposite (face), and one of the three triples listed
above would be a primitive Pythagorean triple. Obviously the triple
$(x,x,z)$ cannot be primitive Pythagorean, since it cannot be Pythagorean
to begin with, since $x^2+x^2=2x^2 \neq z^2$, for any positive integers
$x$ and $z$. So, let us suppose to the contrary that $(x,z,w)$ is a
primitive Pythagorean triple; since equation (\ref{E1}) must also be
satisfied, we must have,

\vspace{.15in}

$
\left\{ \begin{array}{c} 2x^2+z^2=t^2\\
\\
x^2+z^2=w^2 \\
\\
(x,z) = 1 \end{array}\right\} \Leftrightarrow \left\{ \begin{array}{c}
x^2+w^2=t^2\\ \\ x^2+z^2 =w^2\\
\\
(x,z)=1\end{array}\right\}$ \hfill Equations (11)

\vspace{.15in}

\noindent Since $(x,z)=1$, the second equation in (11) describes a
primitive Pythagorean triple, namely $(x,z,w)$; thus, $w$ must be an odd
positive integer.  Moreover, since $(x,z,w)$ is a primitive Pythagorean
triple, we must in fact have $(x,z)=1=(x,w)=(w,z)$; as the condition
$(x,z)=1$ combined with the second equation in (11) implies.  Since
$(x,w)=1$, the first equation in (11) also describes a primitive
Pythagorean triple, namely $(x,w,t)$.  Therefore, by virtue of the fact
that $w$ is odd, the first equation in (11) implies that $x$ must be even
and $t$ odd; and $z$ must be odd.  Altogether,

\vspace{.15in}

\noindent $
\left.\begin{array}{l}
x=2mn,\ w=m^2-n^2,\ t=m^2+n^2\\
\\
x = 2MN,\ z=M^2-N^2,\ w = M^2+N^2,\\
\\
{\rm for\ positive\ integers}\ m,n,M,N,\\
\\
{\rm such\ that\ } m>n,\ M>n,\\
```

\\
(m,n)=1=(M,N),\ {\rm and}\ m+n\equiv 1\equiv M+N\ ({\rm mod}\
2)\end{array}\right\}$ Equations (12)

\vspace{.25in}

\noindent From equations (12) $\Rightarrow \left\{ \begin{array}{l}
mn=MN\\ \\ m^2=n^2 +M^2+N^2 \end{array}\right\}$ \hfill Equations (13)

\vspace{.15in}

\noindent Let $d$ be the greatest common divisor of $m$ and $M;\
(m,M)=d$.

We have,

\vspace{.15in}

$
\left.\begin{array}{r}
m=d\cdot m_1\\
\\
M=d\cdot M_1,\\
\\
{\rm for\ positive\ integers}\ m_1,M_2\ {\rm such}\\
\\
{\rm that}\ (m_1,M_1)=1 \end{array}\right\}$\hfill Equations (14)

\vspace{.15in}

Using the two equations in (14) and the first equation in (13) we obtain
$m_1n=M_1N$; and since $(m_1,M_1)=1$, we see that $m_1$ must be a divisor
of $N$; we set $N=m_1\cdot k$.  So that,

$\left. \begin{array}{r}
N=m_1\cdot k\\
\\
n=M_1 \cdot k,\\
\\
{\rm for\ some\ positive\ integer}\ k.\end{array}\right\}$ \hfill
Equations (15)

\vspace{.15in}

Note that since $(m_1,M_1) = 1,\ k$ is none other that the greatest
common divisor of the integers $n$ and $N$.

Next, by applying the conditions $(m,n)=1=(M,N)$ from (12), together with
(14) and (15) we deduce that,

\setcounter{equation}{15}
\begin{equation}\begin{array}{c}

```latex
\left(M_1\cdot k,d\cdot m_1\right) = 1 = \left(d\cdot M_1,m_1\cdot
k\right) \Rightarrow\\
\\
\Rightarrow \left(M_1,m_1\right) = \left(M_1,d\right) = \left( k,d\right)
= \left(k,M_1\right) = \left(d,M_1\right) = 1
\end{array}
\label{E16}
\end{equation}

Substituting for $m,M,N,n$ from equations (14) and (15) into the second
equation of (13) we obtain,

\begin{equation}
m^2_1 \cdot \left( d^2-k^2\right) = M^2_1 \cdot \left(d^2+k^2\right)
\label{E17}
\end{equation}

On account of $\left(m_1,M_1\right) = 1$ it follows that
$\left(m^2_1,M^2_1\right) = 1$  and consequently (17) implies that
$m^2_1$ must be a divisor of $d^2+k^2$; we set $d^2+k^2=\rho \cdot
m^2_1$; and thus from (17) we now arrive at,

\vspace{.15in}

$ \left. \begin{array}{r}
d^2+k^2 = \rho \cdot m^2_1\\
\\
d^2-k^2 = \rho \cdot M^2_1 \end{array}\right\}$\hfill Equations (18)

\vspace{.15in}

\noindent Furthermore, since by (16) we know that $(d,k)=1$, it follows
that $\left( d^2,k^2)\right) = 1 \Rightarrow \left(d^2+k^2,\ d^2-
k^2\right) = 1$ or $2$, depending on whether the integers $d$ and $k$ are
of different parity or they are both odd.  Clearly, the positive integer
$\rho$ is a common divisor of the integers $d^2+k^2$ and $d^2-k^2$; hence
the only possible values for $\rho$ are $1$ or $2$. If $\rho = 1$, then
according to equations (18), the sum and the difference of two integers
squares, would both be integer squares, which is an impossibility by {\bf
Lemma 1}, which can be also found in \cite{r2}.

Below we show that the case $\rho = 2$ also leads to a contradiction.
Indeed, if $\rho = 2$, the second equation in (18) becomes,

\setcounter{equation}{18}
\begin{equation}
d^2=k^2+2M^2_1 \label{E19}
\end{equation}

The general solution to the diophantine equation $Z^2=Y^2 +2X^2$ is well-
known and can be found in \cite{r3}.

Accordingly then we have,
```

$$
d=D\cdot \left(t^2_1 +2t^2_2\right), \ k = D\cdot \left| t^2_1 - 2t^2_2
\right|,\ M_1 = D\cdot \left(2t_1t_2\right),
$$

\noindent for positive integers $D,t_1,t_2$, such that
$\left(t_1,t_2\right)=1$. Since $(d,k) = 1$, it is obvious that $D=1$,
which gives $d=t^2_1 +2t^2_2,\ k = \left| t^2_1 - 2t^2_2\right|,\ M_1 =
2t_1t_2$. Substituting for $d$ and $k$ in the first equation of (18)
with $\rho =1$ we obtain,

\begin{equation}
\begin{array}{c}
\left(t^2_1+2t^2_2 \right)^2 + \left| t^2_1 - 2t^2_2\right|^2 = 2M^2_1;\\
\\
\left( t^2_1+2t^2_2\right)^2 + \left( t^2_1 - 2t^2_2 \right)^2 = 2M^2_1
\Rightarrow\\
\\
\Rightarrow 2t^4_1+8t^4_2 = 2M^2_1;\\
\\
\left( t^2_1\right)^2 + \left(2t^2_2\right)^2 = M^2_1
\end{array} \label{E20}
\end{equation}

Recall the condition $\left(d,m_1\right)=1$ from (16); this couched with
$d=t^2_1+2t^2_2$ above, shows that $t_1$ must be odd, for if it were even
so would $d$ be, and thus by (20), $m_1$ would also be even, contrary to
$\left(d,M_1\right)=1$. Furthermore, since $t_1$ must be odd and
$\left(t_1,t_2\right) = 1$, we conclude that $\left(t^2_1,2t^2_2\right) =
1$; this shows that equation (20) describes a primitive Pythagorean
triple $\left(t^2_1,2t^2_2,M_1\right)$ with $t^2_1$ odd and thus,

$$
t^2_1 = x^2_1 - x^2_2,\ 2t^2_2=2x_1x_2,\ M_1 = x^2_1 +x^2_2,
$$

\noindent for positive integers $x_1,x_2$ such that $x_1>x_2,\
\left(x_1,x_2\right)=1$ and $x_1+x_2 \equiv 1({\rm mod}\ 2)$. The second
equation rewritten produces $t^2_2=x_1x_2$ and since the integers
$x_1,x_2$ are relatively prime it follows that each must be an integer
square:

$$
x_1 =t^2_3,\ x_2=t^2_4,\ {\rm with}\ \left(t_3,t_4\right) = 1.
$$

\noindent Therefore,

\begin{equation}
t^2_1 = x^2_1 -x^2_2 \Rightarrow t^2_1 = \left(x_1-x_2\right)\left(
x_1+x_2\right) \Rightarrow t^2_1 = \left( t^2_3 -t^2_4\right) \left(t^2_3
+ t^2_4\right). \label{E21}
\end{equation}

Clearly, in view of $\left(x_1,x_2\right)=1$ and $x_1+x_2 \equiv 1({\rm mod}\ 2)$; it follows that $\left(t^2_3,t^2_4\right)=1$ and $t^2_3 + t^2_4 \equiv 1({\rm mod}\ 2)$ which imply $\left( t^2_3 - t^2_4,\ t^2_3 + t^2_4\right) = 1$. Therefore equation (21) shows that both the difference $t^2_3-t^2_4$ and the sum $t^2_3 +t^2_4$ must be integer squares, in violation of {\bf Lemma 1}. \hfill \rule{2mm}{2mm}

\section{ }

\noindent{\bf Lemma 1:}  {\it There do not exist two positive integer squares both of whose sum and difference are also integer squares.}

\vspace{.15in}

\noindent{\bf Proof:}  Let us assume to the contrary that the set $S= \left\{(r,v\left| \right. r,v \in {\mathbb Z}^+,\right.$ $\left.{\rm and\ both}\ r^2-v^2\ {\rm and}\ r^2 + v^2\ {\rm are\ integer\ squares}\right\}$ is nonempty.  Then the set $S_1 = \left\{ n \in {\mathbb Z}^+\left|\right. n=r^2+v^2\ {\rm and\ with}\ (r,v)\in S\right\}$ is a nonempty subset of the set of positive integers or natural numbers; hence it must have a least element, according to the well-ordering principle of the set of natural numbers.  Let $n_0$ be the smallest element of $S_1$; there is a pair $\left(r_0,v_0\right)$ in $S$ such that,

\vspace{.15in}

$ \left. \begin{array}{rl} & r^2_0 - v^2_0 = \beta^2\\
\\
& r^2_0 + v^2_0 = n_0=\gamma ^2,\\
\\
{\rm for\ some} & \beta, \gamma \in {\mathbb Z}^+
\end{array}\right\}$\hfill Equations (22)

\vspace{.15in}

\noindent Let $d$ be the greatest common divisor of $r_0$ and $v_0$; $d = \left( r_0,v_0\right)$. From (22) $\Rightarrow \left( \frac{r_0}{d}\right)^2 - \left( \frac{v_0}{d}\right)^2 = \left( \frac{\beta}{d}\right)^2$ and $\left(\frac{r_0}{d}\right)^2 + \left( \frac{v_0}{d}\right)^2 = \left( \frac{\gamma}{d}\right)^2$, which implies that the pair $\left(\frac{r_0}{d}, \frac{v_0}{d}\right)$ is also a member of $S$.  But $d \geq 1 \Rightarrow \left(\frac{r_0}{d}\right)^2+\left(\frac{v_0}{d} \right)^2 \leq r^2_0 +v^2_0 = n_0$, and since $n_0$ is the least element of $S_1$, it follows that the equal sign must hold in the last inequality, which in turn implies $d = 1 = \left(r_0,v_0\right)$.

Thus both equations in (22) describe primitive Pythagorean triples; the triples $\left(v_0,\beta,r_0\right)$ and $\left( r_0,v_0,\gamma\right)$. Clearly $r_0$ must be odd (from the first triple) and hence (from the second triple) $v_0$ must be even; and both $\beta$ and $\gamma$ must be odd.  If we add the two equations in (22) we obtain,

\setcounter{equation}{22}
\begin{equation}
2r^2_0 = \beta^2+\gamma^2 \Leftrightarrow r^2_0 = \left( \frac{\gamma-\beta}{2}\right)^2 + \left( \frac{\gamma + \beta}{2}\right)^2 \label{E23}
\end{equation}

Note that the two odd numbers $\beta$ and $\gamma$ cannot have a prime divisor in common, for this would imply that the said prime divisor would also divide $r_0$, violating the fact that $\left(v_0,\beta,r_0\right)$ and $\left( r_0,v_0,\gamma \right)$ are primitive Pythagorean triples. Hence $\beta \equiv \gamma \equiv 1({\rm mod}\ 2)$ and $(\beta,\gamma)=1$.  This immediately implies that $\left(\frac{\gamma-\beta}{2},\frac{\gamma + \beta}{2}\right) = 1$; with the two integers $\frac{\gamma - \beta}{2}$ and $\frac{\gamma+\beta}{2}$ having different parity.  Then, according to equation (23), the triple $\left(\frac{\gamma - \beta}{2},\frac{\gamma + \beta}{2}, r_0\right)$ is a primitive Pythagorean one:  Either $\frac{\gamma - \beta}{2} = 2k_1k_2, \frac{\gamma + \beta}{2} = k^2_1-k^2_2, r_0 = k^2_1+k^2_2$; or alternatively $\frac{\gamma - \beta}{2} = k^2_1 - k^2_2, \frac{\gamma +\beta}{2} = 2k_1k_2, r_0 = k^2_1 +k^2_2$, for some natural numbers $k_1,k_2$ such that $k_1 > k_2, \left(k_1,k_2\right) = 1$ and $k_1+k_2 \equiv 1({\rm mod}\ 2)$.  It is clear that in either possibility we must have

$$
\gamma^2-\beta ^2 = 8\left(k_1-k_2\right)\left(k_1+k_2\right) k_1k_2.
$$

Subtracting the first equation from the second in (22) yields,

$$
\gamma^2-\beta^2 = 2v^2_0 \Rightarrow 8 \left(k_1-k_2\right)\left(k_1+k_2\right) k_1k_2 = 2v^2_0 \Rightarrow 4\left(k_1-k_2\right)\left(k_1+k_2\right) k_1k_2 = v^2_0.
$$

In view of $\left(k_1,k_2\right) = 1$ and $k_1+k_2 \equiv 1({\rm mod}\ 2)$, it follows that the four natural nubmers $k_1-k_2,k_1+k_2,k_1,k_2$ must be mutually relatively prime (or coprime); hence, according to the last equation each of them must be an integer square:  $k_1-k_2 = a^2, \ k_1+k_2 = b^2,\ k_1 = c^2,\ k_2 = f^2$, for natural numbers $a,b,c$, and $f$; thus

$$
\left. \begin{array}{rr}
& c^2 - f^2 = a^2\\
\\
{\rm and} & c^2+f^2=b^2 \end{array}\right\}
$$

\noindent which proves that the pair $(c,f)$ is also a member of the set $S$.  However, $c^2+f^2 < r^2_0 + v^2_0 = n_0$, since from the above, it

```
is clear that in fact $c^2+f^2 < v^2_0$.  This is a contradiction since
$n_0$ is the least element of $S_1$.  \hfill \rule{2mm}{2mm}

\vspace{.15in}

\noindent{\bf Remark:}  Even though this proof uses the same idea and
general method, the structure of this proof has a somewhat different
flavor than that found in \cite{r2}.

\section{Closing Remarks and Computations}

\begin{enumerate}
\item[1.]  Note that Result 1 makes use of the diophantine equation
$X^2+Y^2+\delta \cdot Z^2 = \delta \cdot T^2$, for a given odd natural
number $\delta$; one should be able to describe the general solution to
this equation, parametrically and in terms of $\delta$ and its advisors.

\item[2.] With regard to Result 2, below we find the PB which satisfies
Result 2, and which has the smallest diagonal length $t$.  First the
numbers $M_1,N_1,N_2,M_2$ must satisfy equation (1):

$$
M^2_1 +N^2_1 +N^2_2 = M^2_2.
$$

\noindent First assume $M_1$ to be even.  Since the numbers $M_1,N_1$
must also satisfy $M_1>N_1,M_1+N_1 \equiv 1({\rm mod}\ 2)$, and
$\left(M_1,N_1\right) = 1$, then according to the general solution to (1)
we must have, $M_1=2m,\ N_2 = 2\ell,\ N_1 = \frac{\ell^2+m^2-n^2}{n},\
M_2 = \frac{\ell^2+m^2+n^2}{n}$, where $n$ is a natural number which is a
divisor of $\ell ^2 +m^2$.  If we take the first obvious choices of
values for the parameters $\ell $ and $m$; those values being $\ell = 1$
and $m=1$, we then obtain $M_1 = 2, \ N_2 = 2;$ and by also taking $n=1$,
we also obtain $N_1=1$ and $M_2=3$. Indeed,

$$
2^2+2^2+1^2 = 3^2.
$$

Obviously the conditions that $M_1$ and $N_1$ must satisfy, are
satisfied.  Furthermore, $x=2M_1N_1 = 4,\ y=2M_2N_2=12,\ z = M^2_1-N^2_1
= 3$, $w=M^2_1 +N^2_1 = 5$, and $t = M^2_1 +N^2_1 +2N^2_2 = 13$.  check
to see that $x^2+y^2+z^2 = t^2;\ 4^2 +12^2 +3^2 = (13)^2$.  Also
$s=\sqrt{x^2+y^2} = \sqrt{160} = 4\sqrt{10}$ and $u = \sqrt{y^2+z^2} =
\sqrt{153} = 3\sqrt{17}$.  It is now clear that the PB described by
septuple, $(x,y,z,t,s,u,w) = (4,12,3,13, 4\sqrt{10},3\sqrt{17}, 5)$ has a
pair of opposite faces which are primitive PR's; described by the
primitive Pythagorean triple $(x,z,w)=(4,3,5)$.  Also, it is clear that
$t=13$ is the smallest possible value for $t$, among all those PB's
satisfying Result 2, and with $M_1$ even.

Next assume $M_1$ to be odd:  under the conditions $M_1>N_1,\ M_1+N_1
\equiv 1({\rm mod}\ 2),\ (M_1,N_1)=1$, we must have,
```

$$N_1 = 2m,\ N_2 = 2\ell,\ M_1 = \frac{\ell^2 +m^2-n^2}{n},\ M_2 =
\frac{\ell^2+m^2+n^2}{n}.
$$

A search shows that $\ell = m = 2$ and $n=1$ are the values of the parameters $\ell, m$, and $n$ that produce the smallest value of $t$; we obtain $N_1=4,\ N_2=4,\ M_1 = 7,\ M_2 =9$ and consequently $x = 56,\ y = 72,\ z = 33,\ t = 97$.

Thus $t=97$ is the smallest inner diagonal length value, under the assumption that $M_1$ is odd. Hence, $t=13$ is the smallest inner diagonal length value among all PB's generated from Result 2. The PB described by the septuple $(4,12,3,13,4 \sqrt{10},3\sqrt{17}, 5)$ is actually the one with the small $t$ value; $t = 13$, among {\bf all} PB's with a pair of opposite faces which are
primitive PR's; including those PB's that are not generated from Result 2. This should be clear in view of the fact that the triple $(x,z,w)=(4,3,5)$ which describes that pair of PR's, is really the one with the smallest $x,z$ and $w$ values.
\item[3.] Here we ask a question. What is the smallest value of the inner diagonal length $t$, so that both $(x,z,w)$ and $(y,z,u)$ are primitive Pythagorean triples? If we regard the base of such a PB having dimensions $x$ and $y$; then such a PB would have all four lateral faces being primitive PR's and with $t$ being smallest.
\end{enumerate}

\end{document}